\documentclass[12pt, reqno]{amsart}
\usepackage{indentfirst, amssymb, amsmath, amsthm, mathrsfs, setspace, indentfirst, enumerate,  mathrsfs, amsmath, amsthm}

\newtheorem*{theoA}{Theorem A}
\newtheorem*{theoB}{Theorem B}
\newtheorem*{theoC}{Theorem C}

\newcommand{\ol}{\overline}
\newtheorem{theorem}{Theorem}[section]% meant for sectionwise numbers
\newtheorem{lem}{Lemma}[section]
\newcommand{\ee}{\end{equation}}
\newcommand{\beas}{\begin{eqnarray*}}
\newcommand{\eeas}{\end{eqnarray*}}
\newcommand{\bea}{\begin{eqnarray}}
\newcommand{\eea}{\end{eqnarray}}

\numberwithin{equation}{section}
\begin{document}

\title[E\MakeLowercase{ntire Solution of a partial differential equation}]{\LARGE E\Large\MakeLowercase{ntire Solution of a partial differential equation}}

\date{}
\author[J. F. X\MakeLowercase{u}, N. S\MakeLowercase{arkar and} S. M\MakeLowercase{ajumder}]{J\MakeLowercase{unfeng} X\MakeLowercase{u}, N\MakeLowercase{abadwip} S\MakeLowercase{arkar$^*$ and} S\MakeLowercase{ujoy} M\MakeLowercase{ajumder} }
\address{Department of Mathematics, Wuyi University, Jiangmen 529020, Guangdong, People's Republic of China.}
\email{xujunf@gmail.com}
\address{Department of Mathematics, Raiganj University, Raiganj, West Bengal-733134, India.}
\email{naba.iitbmath@gmail.com}
\address{Department of Mathematics, Raiganj University, Raiganj, West Bengal-733134, India.}
\email{sm05math@gmail.com}

\renewcommand{\thefootnote}{}
\footnote{2020 \emph{Mathematics Subject Classification}: }
\footnote{*\emph{Corresponding Author}: Nabadwip Sarkar.}
\renewcommand{\thefootnote}{\arabic{footnote}}
\setcounter{footnote}{0}

\begin{abstract}
In this paper, using Nevanlinna's value distribution theory of meromorphic functions in several complex variables, we study for the existence of entire solutions $f$ in $\mathbb{C}^2$ of the following partial differential equation
\[a_1\left(\frac{\partial f(z_1,z_2)}{\partial z_1}\right)^n+a_2f^n(z_1,z_2)=p_1e^{r(z_1,z_2)}+p_2e^{s(z_1,z_2)},\]
where $n$ is a positive integer such that $n\geq 3$, $a_1,a_2,p_1,p_2$ are non-zero constants and $r(z_1,z_2), s(z_1,z_2)$ are arbitrary polynomials in $\mathbb{C}^2$.
\end{abstract}

\thanks{Typeset by \AmS -\LaTeX}
\maketitle

\section{{\bf Introduction}}
A \textit{Fermat-type functional equation} is a functional equation that generalizes the classical Fermat equation from number theory---namely,
\[
x^{n} + y^{n} = z^{n},
\]
where \(x, y, z\) are typically integers or complex numbers and \(n\) is a positive integer.

When this idea is extended to functions instead of numbers, we obtain \textit{functional equations of Fermat type}, which relate the \(n\)-th powers (or other nonlinear combinations) of functions.

An equation of the form
\bea\label{A}f^n(z)+g^n(z)=1\eea
is a classical example of a \textit{Fermat-type functional equation}. For \( n \geq 2 \), the entire and meromorphic solutions of the equation (\ref{A}) have been thoroughly investigated by Montel \cite{bib4}, Baker \cite{bib0} and Gross \cite{bib1}-\cite{bib3}. To aid the reader, we briefly summarize the established results below.

\begin{theoA} The solutions $f$ and $g$ of the functional equation (\ref{A}) are characterized as follows:
\begin{enumerate}
    \item[(i)] If $ n = 2$,n then the entire solutions are given by $ f = \cos(h) $ and $ g = \sin(h)$, where $ h $ is an entire function and   
    the meromorphic solutions are $f = \frac{1 - \beta^2}{1 + \beta^2}$, $g = \frac{2\beta}{1 + \beta^2}$,
    where $ \beta $ is a non-constant meromorphic function.

    \item[(ii)] If $ n > 2 $, then there are no non-constant entire solutions.

    \item[(iii)] If $ n = 3 $, then the meromorphic solutions are 
 \[f = \frac{1 + \frac{\wp(h(z))}{\sqrt 3}}{2 \wp(h(z))}\;\;\text{and}\;\;g = \frac{1 -\frac{\wp(h(z))}{\sqrt 3}}{2 \wp(h(z))} \, \eta,\]
 where $h$ is a non-constant entire function, $ \eta^3 = 1 $, and $ \wp $ denotes the Weierstrass $\wp$-function satisfying  
    $(\wp')^2 = 4\wp^3 - 1$  
    under appropriate periods.

    \item[(iv)] If $ n > 3$, then there are no non-constant meromorphic solutions.
\end{enumerate}
\end{theoA}

In 1970, Yang \cite{Y1} investigated the following Fermat-type equation
\bea\label{ne2} f^{n}(z)+g^{m}(z)=1\eea
and obtained that equation (\ref{ne2}) has no non-constant entire solutions when $\frac{1}{n}+\frac{1}{m}<1$ (see the proof of Theorem 1). Therefore it is clear that the equation (\ref{ne2}) has no non-constant entire solutions when $n>2$ and $m> 2$. 

In 2004, Yang and Li \cite[Theorem 1]{YL1} considered the Fermat-type differential equation
\bea\label{yl1} f^2(z)+(f^{(1)}(z))^2=1\eea
and proved that if $f(z)$ is a transcendental meromorphic solution of (\ref{yl1}), then $f(z)=\frac{1}{2}\left(Pe^{-\iota z}+\frac{1}{P}e^{\iota z}\right)$, where $P$ is a nonzero constant. If we take $e^A=P$, then solutions of (\ref{yl1}) also can be written as $f(z)=\sin (z+A\iota+\frac{\pi}{2})$.

\medskip
Recently, Fermat-type difference and Fermat-type partial differential-difference equations are the current interests among researchers by utilizing Nevanlinna value distribution theory in several complex variables. We refer the reader to \cite{AA1}-\cite{AM1}, \cite{Chang2012}, \cite{gao2024}, \cite{hal2023}, \cite{hal2024}, \cite{Lu2020}, \cite{MSP}, \cite{MA}, \cite{xu2020b}-\cite{Xu2025}, \cite{zheng2022} as well as the references therein.

\smallskip
In 2020, Xu and Cao \cite{xu2020b} explored transcendental entire functions of finite order that satisfy the following partial differential equation in the variables $(z_1, z_2)$:
\begin{equation}\label{eq1.2}
\left(\frac{\partial f(z_1,z_2)}{\partial z_1}\right)^2 + f^2(z_1, z_2) = 1,
\end{equation}
and obtained the following result.
\begin{theoB}\cite[Corollary 1.4]{xu2020b} The transcendental entire solutions with finite order of the \textit{Fermat-type partial differential equation} (\ref{eq1.2})
have the form $f(z_1, z_2) = \sin(z_1 + g(z_2))$, where $g(z_2)$ is a polynomial in the single variable $z_2$.
\end{theoB}

In 2024, Gao and Liu \cite{gao2024} investigated meromorphic solutions of any order for the following partial differential equation defined in $(z_1, z_2)$:

\begin{equation}\label{eq1.3}
\left( \frac{\partial f(z_1, z_2)}{\partial z_1} \right)^2 + f^2(z_1, z_2) = e^{2r(z_1, z_2)},
\end{equation}
where $r(z_1, z_2)$ is an arbitrary entire function in $\mathbb{C}^2$ and obtained the following result.

\begin{theoC}\cite[Theorem 2.3]{gao2024} Let $r(z)=r(z_1, z_2)$ be a non-constant entire function in $\mathbb{C}^2$. If the equation (\ref{eq1.3})
admits a transcendental meromorphic solution $f(z)$, then one of the following two cases happens,
\begin{itemize}
  \item[(a)] $r(z_1, z_2) = z_1 \cot d + H(z_2)$, then $f(z_1, z_2) = \sin d e^{ z_1 \cot d + H(z_2)},$
  where $H(z_2)$ is an entire function of $z_2$, and $d$ is a complex constant such that $d \neq \frac{k\pi}{2}$ for any $k \in \mathbb{Z}$;
  
  \item[(b)] $r(z_1, z_2) = H(z_2)$ and $f(z_1, z_2) = e^{H(z_2)} \sin(z_1 + T(z_2))$ or $f(z_1, z_2) = \pm e^{H(z_2)},$
  where $H(z_2)$ and $T(z_2)$ are entire functions of $z_2$.
\end{itemize}
\end{theoC}

\section{{\bf Main results}}
A natural question is how to find the entire solutions of the partial differential equation of the form
%Motivated by Theorem B-C, we consider the following partial differential equation:
\begin{equation}\label{eq1.4}
a_1\left( \frac{\partial f(z_1, z_2)}{\partial z_1} \right)^n + a_2f^n(z_1, z_2) = p_1e^{r(z_1, z_2)}+p_2e^{s(z_1, z_2)},
\end{equation}
where $a_1,a_2,p_1,p_2$ are non-zero constants, $n$ is a positive integer and $r(z_1, z_2)$, $s(z_1, z_2)$ are arbitrary polynomials in $\mathbb{C}^2$.
Regarding this question, for $n\geq 5$, we obtain the following result.

\begin{theorem}\label{T3} Let $f(z_1,z_2)$ be an entire solution of the partial differential equation (\ref{eq1.4}). If $n\geq 5$, then $r(z_1,z_2)-s(z_1,z_2)$ reduces to a constant and $f(z_1,z_2)$ is of the form $f(z_1,z_2)=c_1 e^{\frac{s(z_1, z_2)}{n}}$, where $s(z_1,z_2)=az_1+g(z_2)$ and $g(z_2)$ is a polynomial in $z_2$.
\end{theorem}

Therefore a natural question arises: can we find the entire solutions to the equation (\ref{eq1.4}) for $n<5$? However, when $3\leq n\leq 4$ and $r(z_1, z_2)$, $s(z_1, z_2)$ are arbitrary polynomials, obtaining a solution to the equation (\ref{eq1.4}) becomes quite challenging. Therefore, we seek a partial solution by considering
\[r(z_1, z_2) = \lambda_1 z_1 + \gamma_1 z_2\;\;\text{and}\;\;s(z_1, z_2) = \lambda_2 z_1 + \gamma_2 z_2,\]
and obtain the following result.

\begin{theorem}\label{T4}  Let $f$ be an entire solution of Fermat type partial differential equation 
\begin{equation}\label{g01}
a_1\left( \frac{\partial f(z_1, z_2)}{\partial z_1} \right)^n + a_2f^n(z_1, z_2) = p_1e^{\lambda_1 z_1 + \gamma_1 z_2}+p_2e^{\lambda_2 z_1 + \gamma_2 z_2},
\end{equation}
where $n$ is a positive integer such that $3\leq n\leq 4$ and $\lambda_1\neq \lambda_2$. Then $\lambda_2=-3\lambda_1$ or $\lambda_2=\frac{1}{2}(1+\sqrt{3}\iota)\lambda_1$ or $\lambda_2=\frac{1}{2}(1-\sqrt{3}\iota)\lambda_1$ or $\lambda_1+\lambda_2=0$. Moreover, we have
\begin{enumerate}
\item[(1)] for $n=3$ 
\begin{enumerate} 
\item[(i)]  if $\lambda_2=-3\lambda_1$, then  
\[f(z)=c_1(z_2)e^{\lambda_1 z_1}+c_2(z_2)e^{-\lambda_1 z_1},\]
where $c_1(z_2)$ and $c_2(z_2)$ are entire functions such that $a_2c_1^2(z_2)c_2(z_2)=\frac{p_1}{6}e^{\gamma_1 z_2}$, $a_2c_2^3(z)=\frac{p_2}{2}e^{\gamma_2 z_2}$ and $a_1=-\frac{a_2}{\lambda_1^3}$;

\smallskip
\item[(ii)] if $\lambda_2=\frac{1}{2}(1+\sqrt{3}\iota)\lambda_1$, then 
\[f(z_1, z_2) =c_1(z_2)e^{\frac{i\sqrt{3} \lambda_1}{3} z_1}+c_2(z_2)e^{\frac{(3 - i\sqrt{3})\lambda_1}{6} z_1},\]
where $c_1(z_2)$ and $c_2(z_2)$ are entire functions such that $a_2c_1(z_2)c_2^2(z_2)=\frac{p_1 (3+\sqrt 3 i)}{18} e^{\gamma_1 z_2}$, $a_2c_1^2(z_2)c_2(z_2)=\frac{p_2(3-\sqrt 3 i)}{18} e^{\gamma_2 z_2}$ and $a_1 = -\frac{3\sqrt{3}a_2\, i}{\lambda_1^3}$;

\smallskip
\item[(iii)] if $\lambda_2=\frac{1}{2}(1-\sqrt 3 \iota)\lambda_1$, then 
\[f(z_1, z_2)=c_1(z_2)e^{-\frac{i\sqrt{3} \lambda_1}{3} z_1}+c_2(z_2)e^{\frac{(3 + i\sqrt{3})\lambda_1}{6} z_1}\]
where $c_1(z_2)$ and $c_2(z_2)$ are entire functions such that $a_2c_1(z_2)c_2^2(z_2)=\frac{p_1 (3-\sqrt 3 i)}{18} e^{\gamma_1 z_2}$, $a_2c_1^2(z_2)c_2(z_2)=\frac{p_2(3+\sqrt 3 i)}{18} e^{\gamma_2 z_2}$ and $a_1=\frac{3\sqrt{3}a_2\, i}{\lambda_1^3}$,;
\end{enumerate}
\item[(2)] for $n=4$, we have $\lambda_1+\lambda_2=0$ and 
\[f(z_1, z_2) = c_1(z_2)e^{\frac{\lambda_1}{2}z_1}+ c_2(z_2)e^{-\frac{\lambda_1}{2}z_1},\]
where $c_1(z_2)$ and $c_2(z_2)$ are entire functions such that $8a_2c_1^3(z_2)c_2(z)=p_1e^{\gamma_1 z_2}$ and $8a_2c_1(z_2)c_2^3(z)=p_2e^{\gamma_2 z_2}$.
\end{enumerate}

\end{theorem}
%\begin{defi}  Let $f$ and $g$ be meromorphic function in $\mathbb{C}^2$. Then (1) $f=a\Rightarrow g=b$ is defined as zero (s) $f-a$ is/are zero (s) of $g-b$. (2) Zeros of $f$ are algebraic if zeros of $f$ are roots of a polynomial in $\mathbb{C}^2$. 
%\end{defi}

\section {{\bf Notations and Lemmas}}
We define $\mathbb{Z}_+=\mathbb{Z}[0,+\infty)=\{n\in \mathbb{Z}: 0\leq n<+\infty\}$ and $\mathbb{Z}^+=\mathbb{Z}(0,+\infty)=\{n\in \mathbb{Z}: 0<n<+\infty\}$.
On $\mathbb{C}^m$, we define
\[\partial_{z_i}=\frac{\partial}{\partial z_i},\ldots, \partial_{z_i}^{l_i}=\frac{\partial^{l_i}}{\partial z_i^{l_i}}\;\;\text{and}\;\;\partial^{I}=\frac{\partial^{|I|}}{\partial z_1^{i_1}\cdots \partial z_m^{i_m}}\]
where $l_i\in \mathbb{Z}^+\;(i=1,2,\ldots,m)$ and $I=(i_1,\ldots,i_m)\in\mathbb{Z}^m_+$ is a multi-index such that $|I|=\sum_{j=1}^m i_j$.

\smallskip
We firstly recall some basic notions in several complex variables (see \cite{HLY,WS}).
On $\mathbb{C}^m$, the exterior derivative $d$ splits $d= \partial+ \bar{\partial}$ and twists to $d^c= \frac{\iota}{4\pi}\left(\bar{\partial}- \partial\right)$. Clearly $dd^{c}= \frac{\iota}{2\pi}\partial\bar{\partial}$. A non-negative function $\tau: \mathbb{C}^m\to \mathbb{R}[0,b)\;(0<b\leq \infty)$ of class $\mathbb{C}^{\infty}$ is said to be an exhaustion of $\mathbb{C}^m$ if $\tau^{-1}(K)$ is compact whenever $K$ is. 
An exhaustion $\tau_m$ of $\mathbb{C}^m$ is defined by $\tau_m(z)=||z||^2$. The standard Kaehler metric on $\mathbb{C}^m$ is given by $\upsilon_m=dd^c\tau_m>0$. On $\mathbb{C}^m\backslash \{0\}$, we define $\omega_m=dd^c\log \tau_m\geq 0$ and $\sigma_m=d^c\log \tau_m \wedge \omega_m^{m-1}$. For any $S\subseteq \mathbb{C}^m$, let $S[r]$, $S(r)$ and $S\langle r\rangle$ be the intersection of $S$ with respectively the closed ball, the open ball, the sphere of radius $r>0$ centered at $0\in \mathbb{C}^m$.

\smallskip
Let $f$ be a holomorphic function on $G(\not=\varnothing)$, where $G$ is an open subset of $\mathbb{C}^m$. Then we can write $f(z)=\sum_{i=0}^{\infty}P_i(z-a)$, where the term $P_i(z-a)$ is either identically zero or a homogeneous polynomial of degree $i$. Certainly the zero multiplicity $\mu^0_f(a)$ of $f$ at a point $a\in G$ is defined by $\mu^0_f(a)=\min\{i:P_i(z-a)\not\equiv 0\}$.

\medskip
Let $f$ be a meromorphic function on $G$. Then there exist holomorphic functions $g$ and $h$ such that $hf=g$ on $G$ and $\dim_z h^{-1}(\{0\})\cap g^{-1}(\{0\})\leq m-2$. Therefore the $c$-multiplicity of $f$ is just $\mu^c_f=\mu^0_{g-ch}$ if $c\in\mathbb{C}$ and $\mu^c_f=\mu^0_h$ if $c=\infty$. The function $\mu^c_f: \mathbb{C}^m\to \mathbb{Z}$ is nonnegative and is called the $c$-divisor of $f$. If $f\not\equiv 0$ on each component of $G$, then $\nu=\mu_f=\mu^0_f-\mu^{\infty}_f$ is called the divisor of $f$. We define 
$\text{supp}\; \nu=\text{supp}\;\mu_f=\ol{\{z\in G: \nu(z)\neq 0\}}$.

\smallskip
For $t>0$, the counting function $n_{\nu}$ is defined by
\beas n_{\nu}(t)=t^{-2(m-1)}\int_{A[t]}\nu \upsilon_m^{m-1},\eeas
where $A=\text{supp}\;\nu$. 
The valence function of $\nu$ is defined by 
\[N_{\nu}(r)=N_{\nu}(r,r_0)=\int_{r_0}^r n_{\nu}(t)\frac{dt}{t}\;\;(r\geq r_0).\]

Also we write $N_{\mu_f^a}(r)=N(r,a;f)$ if $a\in\mathbb{C}$ and $N_{\mu_f^a}(r)=N(r,f)$ if $a=\infty$.
For $k\in\mathbb{N}$, define the truncated multiplicity function on $\mathbb{C}^m$ by $\mu_{f,k}^a(z)=\min\{\mu_f^a(z),k\}$,
%\beas \mu_{f)k}^a(z)=\begin{cases}
%\mu_f^a(z), &\text{if $\mu_f^a(z)\leq  k$}\\
%0, &\text{if $\mu_f^a(z)>k$}
%\end{cases},\;\;
%\bar{\mu}_{f)k}^a(z)=\begin{cases}
%1, &\text{if $\mu_f^a(z)\leq  k$}\\
%0, &\text{if $\mu_f^a(z)>k$},
%\end{cases}\eeas
%\beas \mu_{f(k}^a(z)=\begin{cases}
%\mu_f^a(z), &\text{if $\mu_f^a(z)\geq  k$}\\
%0, &\text{if $\mu_f^a(z)<k$}
%\end{cases},\;\;
%\bar{\mu}_{f(k}^a(z)=\begin{cases}
%1, &\text{if $\mu_f^a(z)\geq  k$}\\
%0, &\text{if $\mu_f^a(z)<k$}
%\end{cases}\eeas
and the truncated valence functions
\beas N_{\nu}(t)=\begin{cases}
N_k(t,a;f), &\text{if $\nu=\mu_{f,k}^a$}\\
\ol{N}(t,a;f), &\text{if $\nu=\mu_{f,1}^a$}.
%N_{k)}(t,a;f), &\text{if $\nu=\mu_{f)k}^a$}\\
%\ol{N}_{k)}(t,a;f), &\text{if $\nu=\bar{\mu}_{f)k}^a$}\\
%N_{(k}(t,a;f), &\text{if $\nu=\mu_{f(k}^a$}\\
%\ol{N}_{(k}(t,a;f), &\text{if $\nu=\bar{\mu}_{f(k}^a$}.
\end{cases}\eeas

\smallskip
\medskip
An algebraic subset $X$ of $\mathbb{C}^m$ is defined as a subset
\[X=\left\lbrace z\in\mathbb{C}^m: P_j(z)=0,\;1\leq j\leq l\right\rbrace\]
with finitely many polynomials $P_1(z),\ldots, P_l(z)$.
A divisor $\nu$ on $\mathbb{C}^m$ is said to be algebraic if $\nu$ is the zero divisor of a polynomial. In this case the counting
function $n_{\nu}$ is bounded (see \cite{GK1,ST1}). In this case $N_{\nu}(r)=O(\log r)$.

\medskip
With the help of the positive logarithm function, we define the proximity function of $f$ by
\[m(r, f)=\mathbb{C}^m\langle r; \log^+ | f | \rangle=\int_{\mathbb{C}^m\langle r\rangle} \log^+ |f|\;\sigma_m.\]

The characteristic function of $f$ is defined by $T(r, f)=m(r,f)+N(r,f)$. We define $m(r,a;f)=m(r,f)$ if $a=\infty$ and $m(r,a;f)=m(r,1/(f-a))$ if $a$ is finite complex number. Now if $a\in\mathbb{C}$, then the first main theorem of Nevanlinna theory states that $m(r,a;f)+N(r,a;f)=T(r,f)+O(1)$, where $O(1)$ denotes a bounded function when $r$ is sufficiently large. Clearly $f$ is non-constant, then $T(r, f) \rightarrow \infty$ as $r \rightarrow$ $\infty$. Further $f$ is rational if and only if $T(r,f)=O(\log r)$ (see \cite{GK1}).
On the other hand, $f$ to be transcendental then 
\[\limsup\limits_{r \rightarrow \infty} \frac{T(r, f)}{\log r}=+\infty.\] 
We define the order of $f$ by
\[\rho(f):=\limsup _{r \rightarrow \infty} \frac{\log T(r, f)}{\log r}.\]

We define the linear measure $m(E):=\int_E dt$ and the logarithmic measure $l(E):=\int_{E\cap [1,\infty)} \frac{d t}{t}$.

\smallskip
In this paper, for $a\in\mathbb{C}$, we write $f=a\Rightarrow g=a$, if $\{z\in\mathbb{C}^m: f(z)-a=0\}\subseteq \{z\in\mathbb{C}^m: g(z)-a=0\}$.

\medskip
To prove the Theorem \ref{T3} and \ref{T4}, we need the following lemmas.

\smallskip
First we recall the lemma of logarithmic derivative:

\begin{lem}\label{l2} \cite[Lemma 1.37]{HLY} Let $f$ be a non-constant meromorphic function in $\mathbb{C}^m$ and $I=(\alpha_1,\ldots,\alpha_m)\in \mathbb{Z}^m_+$ be a multi-index. Then for any $\varepsilon>0$, we have
\[\parallel\;m\left(r,\frac{\partial^I(f)}{f}\right)\leq |I|\log^+T(r,f)+|I|(1+\varepsilon)\log^+\log T(r,f)+O(1)\]
holds, where $\parallel$ indicates that the inequality holds only outside a set $E$ such that $l(E)<+\infty$.
\end{lem}

\begin{lem}\label{lem1}\cite[Lemma 3.1]{HLY} Let $f_j\not\equiv 0\;(j=1,2,3)$ be meromorphic functions on $\mathbb{C}^m$ such that $f_1$ is non-constant and $f_1+f_2+f_3\equiv 1$.
If
\beas \parallel \;\;\sideset{}{_{j=1}^3}{\sum}\lbrace N_2(r,a_j;f_j)+2\ol{N}(r,f_j)\rbrace<\lambda T(r,f_1)+O(\log^+(T(r,f_1))\eeas
holds, where $\lambda <1$ is positive number, then either $f_2\equiv 1$ or $f_3\equiv 1$.
\end{lem}

%\begin{lem}\label{lem1}\cite{Hu2003}
%Let $f_j \not\equiv 0$ $(j = 1, 2, 3)$ be meromorphic functions on $\mathbb{C}^n$ such that $f_1$ is not constant, $f_1 + f_2 + f_3 = 1$, and such that
%\[
%\sum_{j=1}^3 \left( \overline{N}_2\left(r, \frac{1}{f_j}\right) + 2\overline{N}(r, f_j) \right) \leq \lambda \sum_{j=1}^3 T(r, f_j) + O(\log^+ T(r, f_j))
%\]
%holds for all $r$ outside possibly a set with finite logarithmic measure, where $\lambda < 1$ is a positive number. Then, either $f_2 \equiv 1$ or $f_3 \equiv 1$.
%\end{lem}
\begin{lem}\label{L.2} \cite[Lemma 1.2]{Hu1996} Let $f$ be a non-constant meromorphic function in $\mathbb{C}^m$ and let $a_1,a_2,\ldots,a_q$ be different points in $\mathbb{C}\cup \{\infty\}$. Then
\beas \parallel (q-2)T(r,f)\leq \sideset{}{_{j=1}^{q}}{\sum} \ol N(r,a_j;f)+O(\log (rT(r,f))).\eeas
%where 
%$\parallel$ indicates that the inequality holds only outside a set of finite measure on $\mathbb{R}^+$ and
%$\sideset{}{_{j=1}^{q}}{\sum} N\big(r,\frac{1}{f-a_j}\big)\leq N_{\text{Ram}}(r,f)+\sideset{}{_{j=1}^{q}}{\sum}\ol N\big(r,\frac{1}{f-a_j}\big)$.
\end{lem}

\begin{lem}\label{L.3} \cite[Theorem 1.26]{Hu1996} Let $f$ be non-constant meromorphic function in $\mathbb{C}^m$. Assume that 
$R(z, w)=\frac{A(z, w)}{B(z, w)}$. Then
\beas T\left(r, R_f\right)=\max \{p, q\} T(r, f)+O\Big(\sideset{}{_{j=0}^p}{\sum} T(r, a_j)+\sideset{}{_{j=0}^q}{\sum}T(r, b_j)\Big),\eeas
where $R_f(z)=R(z, f(z))$ and two coprime polynomials $A(z, w)$ and $B(z,w)$ are given
respectively as follows: $A(z,w)=\sum_{j=0}^p a_j(z)w^j$ and $B(z,w)=\sum_{j=0}^q b_j(z)w^j$.
\end{lem}

%\begin{lem}\label{l2}\cite{ye1995} Let \( f \) be a transcendental meromorphic function on \( \mathbb{C}^m \), and let  
%\( I = \{ i_1, i_2, \ldots, i_m \} \) be a multi-index. Then,
%\[
%m\left(r, \frac{\partial^I f}{f} \right) 
%= \int_{S_m(r)} \log^+ \left| \frac{\partial^I f(z)}{f(z)} \right| \sigma_m(z)
%= O\left( \log r \cdot T_f(r) \right)
%\]
%holds for all sufficiently large \( r \), outside a set of finite Lebesgue measure.
%\end{lem}

\section{{\bf Proof of Theorems}}
\subsection{{\bf Proof of Theorem \ref{T3}}}
We prove the theorem by considering the following two cases.\par

\medskip
{\bf Case 1.} Let $r(z_1,z_2)-s(z_1,z_2)$ is a non-constant polynomials. Set 
\[G_1(z_1,z_2) = \frac{a_1}{p_2} \left(f(z_1,z_2) e^{-\frac{s(z_1,z_2)}{n}}\right)^n,\] 
\[G_2(z_1,z_2) = \frac{a_2}{p_2} \left( \frac{\partial f(z_1,z_2)}{\partial z_1} e^{-\frac{s(z_1,z_2) }{n}} \right)^n\]
and
\[G_3(z_1,z_2)= -\frac{p_1}{p_2} e^{(r(z_1,z_2) -s(z_1,z_2)}.\]

Then, from (\ref{eq1.4}), we have
\begin{equation}\label{s2}
G_1 + G_2 + G_3 = 1.
\end{equation}

Note that
\begin{equation}\label{s3}
\begin{cases}
T(r, G_1) =  T\left(r, \frac{a_1}{p_2} \left(f(z_1,z_2)e^{-\frac{s(z_1,z_2)}{n}}\right)^n \right),\\
T(r, G_2) =  T\left(r, \frac{a_2}{p_2}\left( \frac{\partial f(z_1,z_2)}{\partial z_1} e^{-\frac{s(z_1,z_2) }{n}} \right)^n\right),\\
T(r, G_3) = T(r, e^{(r(z_1,z_2) -s(z_1,z_2)}).
\end{cases}
\end{equation}

Also, we have
\begin{equation}\label{s4}
\begin{cases}
N_2(r, 0; G_1) + \overline{N}(r, G_1) \leq 2 \overline{N}\left(r, 0; f(z_1,z_2) e^{-\frac{s(z_1,z_2)}{n}} \right),\\
N_2(r, 0; G_2) + \overline{N}(r, G_2) \leq 2 \overline{N}\left(r, 0; \frac{\partial f(z_1,z_2)}{\partial z_1} e^{-\frac{s(z_1,z_2)}{n}} \right),\\
N_2(r, 0; G_3) + \overline{N}(r, G_3) \leq 2 \overline{N}\left(r, 0; e^{(r(z_1,z_2) -s(z_1,z_2)} \right) = 0.
\end{cases}
\end{equation}

Now, from (\ref{s3}) and (\ref{s4}), we obtain
\beas
&&\sum_{j=1}^{3} N_2(r, 0; G_j) + \sum_{j=1}^{3} \overline{N}(r, G_j)\\ 
&\leq&  2\left( T\left(r, f(z_1,z_2) e^{-\frac{s(z_1,z_2)}{n}} \right) + T\left(r, \frac{\partial f(z_1,z_2)}{\partial z_1} e^{-\frac{s(z_1,z_2)}{n}} \right) \right)\\
&\leq& 2 \left( \frac{1}{n} + \frac{1}{n} \right) T_1(r)
\\&=& \frac{4}{n} T_1(r),
\eeas
where $T_1(r) = \max \left\{ T(r, G_1), T(r, G_2), T(r, G_3) \right\}$. Since $n \geq 5$, using Lemma \ref{lem1}, we conclude that either $G_1 = 1$ or $G_2 = 1$.

\smallskip
First suppose that $G_1 = 1$. Then
\[f(z_1,z_2) = \left(\frac{p_2}{a_1}\right)^{1/n} e^{\frac{s(z_1,z_2)}{n}}\]
and hence
\[\frac{\partial f(z_1,z_2)}{\partial z_1} = \left(\frac{p_2}{a_1}\right)^{1/n} \left(\frac{1}{n} \frac{\partial s(z_1,z_2)}{\partial z_1} \right) e^{\frac{s(z_1,z_2)}{n}}\]
and so from (\ref{eq1.4}), we get 
\[a_2 p_2 \left( \frac{s(z_1,z_2)}{n} \right)^n e^{s(z_1,z_2) - r(z_1,z_2)} = a_1p_1,\]
which leads to a contradiction.\par

\smallskip
Next, suppose that $G_2 = 1$. Then
\begin{equation}\label{s5}
\frac{\partial f(z_1,z_2)}{\partial z_1} = \left( \frac{p_2}{a_2} \right)^{1/n} e^{\frac{s(z_1,z_2) }{n}}.
\end{equation}

Also from (\ref{eq1.4}), we have 
\[f(z_1,z_2) = a\left(\frac{p_1}{a_1}\right)^{1/n} e^{\frac{r(z_1,z_2) }{n}},\]
where $a^n=1$ and so
\begin{equation}\label{s6}
\frac{\partial f(z_1,z_2)}{\partial z_1} = a\left(\frac{p_1}{a_1}\right)^{1/n} \left( \frac{1}{n}\frac{\partial r(z_1,z_2)}{\partial z_1} \right) e^{\frac{r(z_1,z_2)}{n}}.
\end{equation}

Now comparing (\ref{s5}) and (\ref{s6}), we obtain
\[a_2 p_1 \left( \frac{1}{n}\frac{\partial r(z_1,z_2)}{\partial z_1}  \right)^n e^{r(z_1,z_2) - s(z_1,z_2)} = a_1p_2,\]
which leads to a contradiction.\par

\medskip
{\bf Case 2.} Let $r(z_1,z_2)-s(z_1,z_2)$ is constant. Suppose $r(z_1,z_2)-s(z_1,z_2)=c$, where $c\in\mathbb{C}$.
Now equation (\ref{eq1.4}) can be written as follows
\bea\label{cc1} \alpha F^n(z_1,z_2)+\beta G^n(z_1,z_2)=1,\eea
where 
\[\alpha=\frac{a_2}{p_1e^{c}+p_2},\;\;\beta=\frac{a_1}{p_1e^{c}+p_2},\]
\[F(z_1,z_2)=\frac{f(z_1,z_2)}{e^{\frac{s(z_1,z_2)}{n}}}\;\;\text{and}\;\;G(z_1,z_2)=\frac{\frac{\partial f(z_1,z_2)}{\partial z_1}}{e^{\frac{s(z_1,z_2)}{n}}}.\]

Applying Lemma \ref{L.3} to equation (\ref{cc1}), we obtain
\begin{equation}\label{sdd.1}
\parallel\; n T(r,F) + o(T(r,F)) = n T(r,G) + o(T(r,G)).
\end{equation}

Define the function
\begin{equation}\label{sdd.2}
h(z) = \frac{\alpha F^{n}(z_1,z_2) - 1}{\alpha F^{n}(z_1,z_2)}.
\end{equation}

Clearly, \( h \) is a non-constant meromorphic function defined on \( \mathbb{C}^2 \). Again using Lemma \ref{L.3} on (\ref{sdd.2}), we find that
\[\parallel\; T(r,h) + o(T(r,h)) = n T(r,F) + o(T(r,F)).\]

Now from equations (\ref{cc1}) and (\ref{sdd.2}), it follows that
\[
\parallel\; \ol{N}(r,h) \leq \ol{N}(r,0, F^{n}) = \ol{N}(r,0,F) + o(T(r,F)),
\]
\beas
\parallel\; \ol{N}(r,0,h) &=& \ol{N}(r,1/\alpha, F^{n}) \leq \ol{N}(r,0, G^{n})= \ol{N}(r,0,G) + o(T(r,G)),
\eeas
\[
\parallel\; \ol{N}(r,1,h) \leq \ol{N}(r, F^{n}) =  o(T(r,F)).
\]

Thus, by Lemma \ref{L.2}, we derive the inequality
\begin{equation}\label{sdd.3}
\begin{aligned}
\parallel\; n T(r,F) &= T(r,h) + o(T(r,h)) \\
&\leq \ol{N}(r,h) + \ol{N}(r,0,h) + \ol{N}(r,1,h) + o(T(r,h)) \\
&\leq \ol{N}(r,0,F) + \ol{N}(r,0,G)  + o(T(r,F)) + o(T(r,G)).
\end{aligned}
\end{equation}

Now, by the first main theorem and using (\ref{sdd.1}) to (\ref{sdd.3}), we get
\[\parallel\; \left(n - 2 \right) T(r,F) \leq o(T(r,F)),\]
which shows that $F$ is constant, say $c_1$. Since $F(z_1,z_2)=f(z_1,z_2)e^{-\frac{s(z_1,z_2)}{n}}$, it follows that
\bea\label{xsm1}f(z_1, z_2) = c_1 e^{\frac{s(z_1, z_2)}{n}}.\eea
 
Again from (\ref{sdd.1}), we can deduce that \( G \) is also a constant, say \( c_2 \). Since $G(z_1,z_2)=\frac{\partial f(z_1,z_2)}{\partial z_1}e^{-\frac{s(z_1,z_2)}{n}}$, it follows that
\bea\label{xsm2}\frac{\partial f(z_1, z_2)}{\partial z_1} = c_2 e^{\frac{s(z_1, z_2)}{n}}.\eea

Now from (\ref{xsm1}) and (\ref{xsm2}), we get 
\[\frac{\partial s(z_1, z_2)}{\partial z_1} = n \frac{c_2}{c_1},\]
i.e., $s(z_1,z_2)$ is of the form \( s(z_1, z_2) = a z_1 + g(z_2) \), where \( g(z_2) \) is a polynomial in \( z_2 \).
Finally we have $f(z)=c_1 e^{\frac{s(z_1, z_2)}{n}}$, where $s(z_1,z_2)=az_1+g(z_2)$ such that $g(z_2)$ is a polynomial in $z_2$ and $a,c_1$ are non-zero constants. 

Hence the proof is complete.

\subsection{{\bf Proof of Theorem \ref{T4}}}
Let \( f \) be an entire solution of the equation (\ref{g01}). Clearly, \( f \) is transcendental. We rewrite the (\ref{g01}) as follows
\begin{equation}\label{g1}
\hat a\left( \frac{\partial f(z_1, z_2)}{\partial z_1} \right)^n + f^n(z_1, z_2) = \hat p_1e^{\lambda_1 z_1 + \gamma_1 z_2}+\hat p_2e^{\lambda_2 z_1 + \gamma_2 z_2},
\end{equation}
where $\hat a=\frac{a_1}{a_2}$, $\hat p_1=\frac{p_1}{a_2}$ and $\hat p_2=\frac{p_2}{a_2}$.
Now differentiating (\ref{g1}) partially with respect to $z_1$, we get
\begin{equation} \label{g2}
n\left(f^{n-1} \frac{\partial f}{\partial z_1} + \hat a \left(\frac{\partial f}{\partial z_1}\right)^{n-1} \frac{\partial^2 f}{\partial z_1^2} \right) = \hat p_1 \lambda_1 e^{\lambda_1 z_1+\gamma_1z_2} + \hat p_2 \lambda_2 e^{\lambda_2 z_1+\gamma_2z_2}.
\end{equation}

Eliminating \( e^{\lambda_1 z_1+\gamma_1z_2} \) from (\ref{g1}) and (\ref{g2}), we get
\bea \label{g3}&& n\left(f^{n-1} \frac{\partial f}{\partial z_1} + \hat a \left(\frac{\partial f}{\partial z_1}\right)^{n-1} \frac{\partial^2 f)}{\partial z_1^2} \right) - \lambda_1\left(f^n + \hat a \left(\frac{\partial f}{\partial z_1}\right)^n\right) \nonumber\\
&&= \hat p_2 (\lambda_2 - \lambda_1) e^{\lambda_2 z_1+\gamma_2z_2}.\eea

Differentiating partially (\ref{g3}) with respect to $z_1$, we get
\bea \label{g4}
&& n \left((n-1) f^{n-2} \left(\frac{\partial f}{\partial z_1}\right)^2 + f^{n-1} \frac{\partial^2 f}{\partial z_1^2} \right)\nonumber\\
&& +\hat  a n \left((n-1)\left(\frac{\partial f}{\partial z_1}\right)^{n-2} \left(\frac{\partial^2 f}{\partial z_1^2}\right)^2 + \left(\frac{\partial f}{\partial z_1}\right)^{n-1} \frac{\partial^3 f}{\partial z_1^3} \right) \nonumber\\
&& - n \lambda_1 \left(f^{n-1} \frac{\partial f}{\partial z_1} + \hat a \left(\frac{\partial f}{\partial z_1}\right)^{n-1} \frac{\partial^2 f}{\partial z_1^2} \right)= \hat p_2 \lambda_2 (\lambda_2 - \lambda_1) e^{\lambda_2 z_1+\gamma_2z_2}.
\eea

Now eliminating \( e^{\lambda_2 z_1+\gamma_2z_2} \) from (\ref{g3}) and (\ref{g4}), we get
\bea \label{g6}
&&\left(n \frac{\partial^2 f}{\partial z_1^2} + \lambda_1 \lambda_2 f \right) f^{n-1}\nonumber \\
&& = \frac{\partial f}{\partial z_1} \Bigg[ 
n (\lambda_1 + \lambda_2) \left(f^{n-1} + \hat a \left(\frac{\partial f}{\partial z_1}\right)^{n-2} \frac{\partial^2 f}{\partial z_1^2} \right) \nonumber\\
&& \quad - n \left( (n-1) f^{n-2} \frac{\partial f}{\partial z_1} + (n-1) \hat a \left(\frac{\partial f}{\partial z_1}\right)^{n-3} \left(\frac{\partial^2 f}{\partial z_1^2}\right)^2 + \hat a \left(\frac{\partial f}{\partial z_1}\right)^{n-2} \frac{\partial^3 f}{\partial z_1^3} \right)\nonumber \\
&& \quad - \lambda_1 \lambda_2 \hat a \left(\frac{\partial f}{\partial z_1}\right)^{n-1}
\Bigg].
\eea

We claim that \( \rho(f) < +\infty \). Suppose, on the contrary, that \( f \) is a solution of infinite order to the equation (\ref{g1}).
Then form every $M>0$ there exists $r_0$ such that 
\bea\label{tt1} T(r,f)>r^M\;\;\text{when}\;r>r_0.\eea 

Note that 
\[T(r,e^{\lambda_1 z_1 + \gamma_1 z_2})=m(r,e^{\lambda_1 z_1 + \gamma_1 z_2})=O(r)\] 
and 
\[T(r,e^{\lambda_2 z_1 + \gamma_2 z_2})=m(r,e^{\lambda_2 z_1 + \gamma_2 z_2})=O(r).\]

If we take 
\[\alpha(z_1, z_2) = \hat p_1 e^{\lambda_1 z_1 + \gamma_1 z_2} + \hat p_2 e^{\lambda_2 z_1 + \gamma_2 z_2},\]
then by Lemma \ref{L.3}, we get 
\[T(r,\alpha(z_1, z_2))\leq T(r,e^{\lambda_1 z_1 + \gamma_1 z_2})+T(r,e^{\lambda_2 z_1 + \gamma_2 z_2})+O(\log r)\]
and so
\bea\label{tt2} T(r,\alpha(z_1, z_2))=O(r).\eea

Let us take $M>1$. Then from \ref{tt1}) and (\ref{tt2}), we get
\[\frac{T(r,\alpha)}{T(r,f)}<\frac{O(r)}{r^M}\to 0,\]
as $r\to \infty$. Therefore we conclude that $\alpha$ is small function of $f$, i.e., $T(r,\alpha)=o(T(r,f))$.

Let \( w_i \; (i = 1, 2, \ldots, n) \) be the roots of the equation $\hat az^n +1 = 0$. Clearly $w_i\neq 0$ for $i=1,2,\ldots,n$. Set
\[g = \frac{\frac{\partial f}{\partial z_1}}{ f}.\]

\medskip
First we suppose that $g$ is a small function of $f$. Then from (\ref{g1}), we have
\[f^n (\hat ag^n + 1)=\alpha\]
and so by Lemma \ref{L.3}, we get a contradiction.

\medskip
Next we suppose that $g$ is not a small function of $f$. Then from (\ref{g1}), we have
\begin{equation} \label{xx4}
\hat a(g - w_1)(g - w_2) \cdots (g - w_n) = \alpha f^{-n}.
\end{equation}

If \( g = w_i \Rightarrow f = 0 \), then obviously \( f = 0 \Rightarrow \frac{\partial f}{\partial z_1} = 0 \), for \( i = 1, 2, \ldots, n \). Clearly, from (\ref{g1}), it follows that \( g = w_i \Rightarrow \alpha = 0 \). 

If \( g = w_i \not\Rightarrow f = 0 \), then from (\ref{xx4}), we again get \( g = w_i \Rightarrow \alpha = 0 \). 

Therefore, in either case, we have \( g = w_i \Rightarrow \alpha = 0 \) for $i=1,2,\ldots,n$.
Consequently using the first main theorem, we have
\begin{equation} \label{h1}
N(r,w_i;g)\leq N(0;\alpha)\leq T(r,\alpha)+O(1)=o(T(r,f)) \quad (i = 1, 2, \ldots, n).
\end{equation}

Now in view of (\ref{h1}) and using Lemma \ref{L.2}, we obtain $T(r,g)=o(T(r, f))$, which contradicts the fact that $g$ is not a small function of $f$.

Finally we deduce that $\rho(f)<+\infty$. Now from equation~(\ref{g6}), we obtain
\bea \label{g6.1}
\lambda_1 \lambda_2 \frac{f}{\frac{\partial f}{\partial z_1}} &=& 
n(\lambda_1 + \lambda_2) \left(1 + \hat a \left( \frac{\frac{\partial f}{\partial z_1}}{f} \right)^{n-2} \frac{\frac{\partial^2 f}{\partial z_1^2}}{f} \right) \nonumber\\
&& - n \left( (n - 1) \frac{\frac{\partial f}{\partial z_1}}{f} 
+ (n - 1) \hat a \left( \frac{\frac{\partial f}{\partial z_1}}{f} \right)^{n - 3} \left( \frac{\frac{\partial^2 f}{\partial z_1^2}}{f} \right)^2 
+ \hat a \left( \frac{\frac{\partial f}{\partial z_1}}{f} \right)^{n - 2} \frac{\frac{\partial^3 f}{\partial z_1^3}}{f} \right) \nonumber\\
&& - \lambda_1 \lambda_2 \hat a \left( \frac{\frac{\partial f}{\partial z_1}}{f} \right)^{n - 1} 
- n \frac{\frac{\partial^2 f}{\partial z_1^2}}{\frac{\partial f}{\partial z_1}}.
\eea

Since $\rho(f)<+\infty$, applying Lemma \ref{l2} to (\ref{g6.1}), we deduce that
\begin{equation} \label{g6.2}
\parallel\;m\left(r, \frac{f}{\frac{\partial f}{\partial z_1}} \right) = O(\log r).
\end{equation}

Now we divide the following two cases.\par

\medskip
{\bf Case 1.} Suppose that \(0 \) is a Picard exceptional value of \( \frac{\partial f}{\partial z_1} \). 
Since \( 0 \) is a Picard exceptional value of \( \frac{\partial f}{\partial z_1} \), it follows that $\frac{f}{\frac{\partial f}{\partial z_1}}$ is an entire function. Now from (\ref{g6.2}), we deduce that 
\[\parallel\;T\left(r,\frac{f}{\frac{\partial f}{\partial z_1}}\right)=O(\log r),\]
which shows that $\zeta=\frac{f}{\frac{\partial f}{\partial z_1}}$ is a non-zero polynomial in $\mathbb{C}^2$. Then from (\ref{g1}), we have
\bea\label{xu1}f^n \left(\frac{\hat a}{\zeta^n} + 1\right) =\hat p_1 e^{\lambda_1 z_1 + \gamma_1 z_2} + \hat p_2 e^{\lambda_2 z_1 + \gamma_2 z_2}=\hat p_2 e^{\lambda_2 z_1 + \gamma_2 z_2}\left(1+\frac{p_1}{p_2}e^{(\lambda_1-\lambda_2)z_1+(\gamma_1-\gamma_2)z_2}\right).\eea

Clearly $\hat a+\zeta^n$ is a non-zero polynomial in $\mathbb{C}^2$. Since the right hand side of (\ref{xu1}) is an entire function in $\mathbb{C}^2$, it follows that if $z_0\in\mathbb{C}^2$ is a zero of $\zeta$ of multiplicity $l$ then $z_0$ must be a zero of $f^n$ of multiplicity at least $nl$. Let 
\[\beta=\frac{p_1}{p_2}e^{(\lambda_1-\lambda_2)z_1+(\gamma_1-\gamma_2)z_2}.\]

Suppose $\zeta\in\mathbb{C}^2$ is a zero of $\beta+1$. Then from (\ref{xu1}), we have the following two possibilities
\begin{enumerate}
\item[(i)] $\beta(\zeta)+1=0\Rightarrow \hat a+\zeta^n(z_1)=0$. In this case, the set of zeros of $\hat a+\zeta^n$ is algebraic;
\item[(ii)] $\beta(\zeta)+1=0\Rightarrow f^n(\zeta)=0$. In this case, all the zeros of $\beta+1$ have multiplicities at least $3$.
\end{enumerate}

\medskip
Now using Lemma \ref{L.2}, we get
\beas \parallel\;T\left(r,\beta\right)&\leq& \ol N\left(r,0;\beta\right)+\ol N(r,\beta)+\ol N\left(r,-1;\beta\right)+o\left(T\left(r,\beta\right)\right)\\&\leq&
\ol N\left(r,0;\hat a+\zeta^n\right)+\frac{1}{3} N\left(r,-1;\beta\right)+o\left(T\left(r,\beta\right)\right)\\&\leq&
\frac{1}{3}N\left(r,-1;\beta\right)+o\left(T\left(r,\beta\right)\right),\nonumber\\&\leq&
\frac{1}{3}T\left(r,\beta\right)+o\left(T\left(r,\beta\right)\right),\nonumber
\eeas
which is impossible.\par

\medskip
{\bf Case 2.} Suppose that $0$ is not Picard exceptional value of $\frac{\partial f}{\partial z_1}$. 

Now we divide the following two sub-cases:\par

\medskip
{\bf Sub-case 2.1.} Let $n =3$. Then from (\ref{g6}), we get
\begin{equation} \label{g7}
\begin{aligned}
\left(3 \frac{\partial^2 f}{\partial z_1^2} + \lambda_1 \lambda_2 f \right) f^2 
&= \frac{\partial f}{\partial z_1} \Bigg[
3(\lambda_1 + \lambda_2)\left( f^2 + \hat a \frac{\partial f}{\partial z_1} \frac{\partial^2 f}{\partial z_1^2} \right) \\
&\quad - 3\left( 2 f \frac{\partial f}{\partial z_1} + 2\hat a \left( \frac{\partial^2 f}{\partial z_1^2} \right)^2 + \hat a \frac{\partial f}{\partial z_1} \frac{\partial^3 f}{\partial z_1^3} \right)
\quad - \lambda_1 \lambda_2 \hat a \left( \frac{\partial f}{\partial z_1} \right)^2 \Bigg].
\end{aligned}
\end{equation}

Let
\begin{equation} \label{g8}
\phi = \frac{3 \frac{\partial^2 f}{\partial z_1^2} + \lambda_1 \lambda_2 f}{\frac{\partial f}{\partial z_1}}.
\end{equation}

Now we consider following two sub-cases:

\smallskip
{\bf Sub-case 2.1.1.} Let \( \phi \equiv 0 \). Then 
\bea\label{xu2} 3 \frac{\partial^2 f}{\partial z_1^2} + \lambda_1 \lambda_2 f \equiv 0\eea
and so from (\ref{g7}), we obtain
\bea\label{xu3}
\left( 3(\lambda_1 + \lambda_2) - \frac{2\hat a}{3} \lambda_1^2 \lambda_2^2 \right) f = \left( (\lambda_1 + \lambda_2) \lambda_1 \lambda_2 \hat a + 6 \right) \frac{\partial f}{\partial z_1}.\eea

\medskip
First we suppose that $3(\lambda_1 + \lambda_2) - \frac{2\hat a}{3} \lambda_1^2 \lambda_2^2\neq 0$. Then from (\ref{xu3}), we have
\[(\lambda_1 + \lambda_2) \lambda_1 \lambda_2 \hat a + 6\neq 0.\]

In this case, from (\ref{xu3}), we see that $\frac{f}{\frac{\partial f}{\partial z_1}}$ is a non-zero constant. Now proceeding in the same way as done in the proof of Case 1, we get a contradiction.

\medskip
Next we suppose that $3(\lambda_1 + \lambda_2) - \frac{2\hat a}{3} \lambda_1^2 \lambda_2^2=0$. Then from (\ref{xu3}), we have
\[(\lambda_1 + \lambda_2) \lambda_1 \lambda_2 \hat a + 6=0.\]

Therefore we have
\begin{equation} \label{kj1}
3(\lambda_1 + \lambda_2) - \frac{2\hat a}{3} \lambda_1^2 \lambda_2^2 = 0 \quad \text{and} \quad (\lambda_1 + \lambda_2)\lambda_1\lambda_2 \hat a + 6 = 0.
\end{equation}

Solving (\ref{kj1}), we get \( \lambda_2 = -3\lambda_1 \) or \( \lambda_1 = -3\lambda_2 \). Without loss of generality, assume \( \lambda_2 = -3\lambda_1 \). Then from (\ref{kj1}), we have \( \hat a = -\frac{1}{\lambda_1^3} \). 

Now, from (\ref{xu2}), we get
\[
\frac{\partial^2 f}{\partial z_1^2} - \lambda_1^2 f = 0,
\]
and hence the general solution is of the form
\[
f(z_1, z_2) = c_1(z_2)e^{\lambda_1z_1}+c_2(z_2)e^{-\lambda_1z_1},
\] 
where $c_1(z_2)$ and $c_2(z_2)$ are entire functions.
Now putting the value of $f$ in (\ref{g1}), we get $c_1^2(z_2)c_2(z_2)=\frac{\hat p_1}{6}e^{\gamma_1 z_2}$ and $c_2^3(z)=\frac{\hat p_2}{2}e^{\gamma_2 z_2}$. 
Finally we have 
\[f(z_1, z_2) = c_1(z_2)e^{\lambda_1z_1}+c_2(z_2)e^{-\lambda_1z_1},\]
where $c_1(z_2)$ and $c_2(z_2)$ are entire functions such that $c_1^2(z_2)c_2(z_2)=\frac{\hat p_1}{6}e^{\gamma_1 z_2}$ and $c_2^3(z)=\frac{\hat p_2}{2}e^{\gamma_2 z_2}$.\par

\medskip
{\bf Sub-case 2.1.2.} Let \( \phi \not\equiv 0 \). Clearly, from (\ref{g3}), we observe that \( f \) and \( \frac{\partial f}{\partial z_1} \) cannot have common zeros. Therefore, from (\ref{g7}), we conclude that if $z_0\in\mathbb{C}^2$ is a zero of $\frac{\partial f}{\partial z_1}$ of multiplicity $l$, then $z_0$ must be a zero of $3\frac{\partial^2 f}{\partial z_1^2} + \lambda_1 \lambda_2 f$ of multiplicity at least $l$.
Consequently \( \phi \) is a non-zero entire function. Now in view of (\ref{g6.2}) and using Lemma \ref{l2} to (\ref{g8}), we obtain
\[\parallel\;m(r, \phi) = O(\log r).\]

Hence, \( \phi \) is a non-zero polynomial in $\mathbb{C}^2$. From equation (\ref{g8}), we obtain:
\begin{equation} \label{g9.1}
\frac{\partial^2 f}{\partial z_1^2} = \frac{1}{3} \left( \phi \frac{\partial f}{\partial z_1} - \lambda_1 \lambda_2 f \right).
\end{equation}

Substituting (\ref{g8}) and (\ref{g9.1}) into (\ref{g7}), we get:
\begin{equation} \label{g10}
\begin{aligned}
&\left( \phi - \left( 3(\lambda_1 + \lambda_2) - \frac{2\hat a}{3} (\lambda_1 \lambda_2)^2 \right) \right) f^2 = f \frac{\partial f}{\partial z_1} \left( -\hat a(\lambda_1 + \lambda_2)\lambda_1 \lambda_2 - 6 + \frac{4\hat a}{3} \phi \lambda_1 \lambda_2 \right) \\
&\quad + \left( \frac{\partial f}{\partial z_1} \right)^2 \left(\hat a(\lambda_1 + \lambda_2)\phi - \frac{2a}{3} \phi^2 -\hat a \lambda_1 \lambda_2 \right) - 3a \frac{\partial f}{\partial z_1} \frac{\partial^3 f}{\partial z_1^3}.
\end{aligned}
\end{equation}

\medskip
First we suppose that
\[\phi - \left( 3(\lambda_1 + \lambda_2) - \frac{2\hat a}{3} (\lambda_1 \lambda_2)^2 \right) \not\equiv 0.\]

Since \( f \) and \( \frac{\partial f}{\partial z_1} \) have no common zeros, it follows from (\ref{g10}) that
\[\frac{\partial f}{\partial z_1} = 0 \Rightarrow \phi - \left( 3(\lambda_1 + \lambda_2) - \frac{2\hat a}{3} (\lambda_1 \lambda_2)^2 \right) = 0.\]

Therefore the set of zeros of \( \frac{\partial f}{\partial z_1} \) is algebraic. Let us define a function 
\beas \varphi=\frac{f}{\frac{\partial f}{\partial z_1}}.\eeas

Clearly $N(r,\varphi)=O(\log r)$. Also in view of (\ref{g6.2}), we see that $\parallel\;m(r, \varphi)=O(\log r)$. Consequently $\varphi$ is a rational function in $\mathbb{C}^2$. Then from (\ref{g1}), we have
\[f^3 \left(\frac{\hat a}{\varphi^3}+1\right)=\hat p_1 e^{\lambda_1 z_1 + \gamma_1 z_2} + \hat p_2 e^{\lambda_2 z_1 + \gamma_2 z_2}=\hat p_2 e^{\lambda_2 z_1 + \gamma_2 z_2}\left(1+\frac{p_1}{p_2}e^{(\lambda_1-\lambda_2)z_1+(\gamma_1-\gamma_2)z_2}\right).\]

Now proceeding in the same way as done in the proof of Case 1, we get a contradiction.

\medskip
Next we suppose that 
\[\phi - \left( 3(\lambda_1 + \lambda_2) -\frac{2\hat a}{3}(\lambda_1 \lambda_2)^2\right)=0,
\quad \text{i.e.,} \quad 
\phi = \frac{1}{3}\left(3(\lambda_1 + \lambda_2) -2\hat a(\lambda_1 \lambda_2)^2\right).\]

Again from (\ref{g10}), we get
\begin{equation} \label{xu4}
\begin{aligned}
&f \left(-\hat a(\lambda_1 + \lambda_2)\lambda_1 \lambda_2 - 6 + \frac{4\hat a}{3} \phi \lambda_1 \lambda_2 \right) 
+ \frac{\partial f}{\partial z_1} \left(\hat a(\lambda_1 + \lambda_2)\phi - \frac{2\hat a}{3} \phi^2 -\hat a\lambda_1 \lambda_2 \right) \\
&\quad -3\hat a \frac{\partial^3 f}{\partial z_1^3}=0.
\end{aligned}
\end{equation}

Differentiating partially (\ref{g9.1}) with respect to $z_1$ and then eliminating $\frac{\partial^3 f}{\partial z_1^3}$ from (\ref{xu4}), we obtain
\[
f \left(-\hat a(\lambda_1 + \lambda_2)\lambda_1 \lambda_2 - 6 + \frac{5\hat a}{3} \phi \lambda_1 \lambda_2 \right)
+ \frac{\partial f}{\partial z_1} \left(\hat a(\lambda_1 + \lambda_2)\phi -\hat a \phi^2 \right) = 0,
\]
which yields
\begin{equation}\label{g11}
-\hat a(\lambda_1 + \lambda_2)\lambda_1 \lambda_2 - 6 + \frac{5\hat a}{3} \phi \lambda_1 \lambda_2 = 0
\quad \text{and} \quad
\hat a(\lambda_1 + \lambda_2)\phi -\hat a \phi^2 = 0.
\end{equation}

Now from (\ref{g11}), we get 
\[\phi = \lambda_1 + \lambda_2\;\;\text{and}\;\;\hat a(\lambda_1 + \lambda_2)\lambda_1 \lambda_2 - 9 = 0.\]

Again from (\ref{g11}), we obtain
\[
3(\lambda_1 + \lambda_2) = \hat a(\lambda_1 \lambda_2)^2.
\]
Therefore, we can easily deduce 
\[
\lambda_2 = \frac{1}{2}(1 + \sqrt{3}\, i)\lambda_1 \quad \text{or} \quad 
\lambda_2 = \frac{1}{2}(1 - \sqrt{3}\, i)\lambda_1.
\]

\smallskip
If $\lambda_2 = \frac{1}{2}(1 + \sqrt{3}\, i)\lambda_1$, then from (\ref{g8}) we get
\[
6 \frac{\partial^2 f}{\partial z_1^2} - (3 + i\sqrt{3})\lambda_1 \frac{\partial f}{\partial z_1} + (1 + i\sqrt{3})\lambda_1^2 f = 0,
\]
and so $f$ is of the form
\[
f(z_1, z_2) =c_1(z_2)e^{\frac{i\sqrt{3} \lambda_1}{3} z_1}+c_2(z_2)e^{\frac{(3 - i\sqrt{3})\lambda_1}{6} z_1}.\]

Now putting the value of $f$ in (\ref{g1}), we get $\hat a = -\frac{3\sqrt{3}\, i}{\lambda_1^3}$, $c_1(z_2)c_2^2(z_2)=\frac{\hat p_1 (3+\sqrt 3 i)}{18} e^{\gamma_1 z_2}$ and $c_1^2(z_2)c_2(z_2)=\frac{\hat p_2(3-\sqrt 3 i)}{18} e^{\gamma_2 z_2}$.

\smallskip
If $\lambda_2 = \frac{1}{2}(1 - \sqrt{3}\, i)\lambda_1$, then from (\ref{g8}) we get
\[
6 \frac{\partial^2 f}{\partial z_1^2} - (3 - i\sqrt{3})\lambda_1 \frac{\partial f}{\partial z_1} + (1 - i\sqrt{3})\lambda_1^2 f = 0,
\]
and so $f$ is of the form
\[
f(z_1, z_2)=c_1(z_2)e^{-\frac{i\sqrt{3} \lambda_1}{3} z_1}+c_2(z_2)e^{\frac{(3 + i\sqrt{3})\lambda_1}{6} z_1}.\]

Now putting the value of $f$ in (\ref{g1}), we get $\hat a = \frac{3\sqrt{3}\, i}{\lambda_1^3}$, $c_1(z_2)c_2^2(z_2)=\frac{\hat p_1 (3-\sqrt 3 i)}{18} e^{\gamma_1 z_2}$ and $c_1^2(z_2)c_2(z_2)=\frac{\hat p_2(3+\sqrt 3 i)}{18} e^{\gamma_2 z_2}$.\par

\medskip
{\bf Case 2.2.} Let $n=4$. Then from (\ref{g6}), we obtain 
\bea\label{g13}
&&\left(4 \frac{\partial^2 f}{\partial z_1^2} + \lambda_1 \lambda_2 f\right) f^3 \nonumber\\
&&= \frac{\partial f}{\partial z_1} \bigg[4(\lambda_1+\lambda_2)\left(f^3 +\hat a \left(\frac{\partial f}{\partial z_1}\right)^2 \frac{\partial^2 f}{\partial z_1^2}\right)\nonumber\\
&&-4\left(3 f^2 \frac{\partial f}{\partial z_1} + 3\hat a \frac{\partial f}{\partial z_1} \left(\frac{\partial^2 f}{\partial z_1^2}\right)^2 +\hat a \left(\frac{\partial f}{\partial z_1}\right)^2 \frac{\partial^3 f}{\partial z_1^3}\right)
- \lambda_1 \lambda_2 \hat a \left(\frac{\partial f}{\partial z_1}\right)^3\bigg].
\eea

Define
\begin{equation}\label{g14}
\psi = \frac{4 \frac{\partial^2 f}{\partial z_1^2} + \lambda_1 \lambda_2 f}{\frac{\partial f}{\partial z_1}}.
\end{equation}

Now we consider following two sub-cases:

\smallskip
{\bf Sub-case 2.2.1.} Suppose that $\psi \equiv 0$. Then from (\ref{g14}), we have 
\bea\label{xu5}
\frac{\partial^2 f}{\partial z_1^2} + \frac{\lambda_1 \lambda_2}{4} f = 0.\eea

\medskip
First we suppose that $\lambda_1 + \lambda_2 \neq 0$. Then from (\ref{g13}), we obtain
\beas 
4(\lambda_1+\lambda_2)f^3 &=& \frac{\partial f}{\partial z_1} \bigg[-4(\lambda_1+\lambda_2)\hat a \frac{\partial f}{\partial z_1} \frac{\partial^2 f}{\partial z_1^2}\\
&&-4\left(3 f^2 + 3\hat a \left(\frac{\partial^2 f}{\partial z_1^2}\right)^2 +\hat a \frac{\partial f}{\partial z_1} \frac{\partial^3 f}{\partial z_1^3}\right)
- \lambda_1 \lambda_2\hat a \left(\frac{\partial f}{\partial z_1}\right)^2 \bigg],
\eeas
which contradicts the fact that $f$ and $\frac{\partial f}{\partial z_1}$ have no common zeros. 

\medskip
Next we suppose that $\lambda_1 + \lambda_2 = 0$. Then from (\ref{xu5}), we have 
\bea\label{xu6}\frac{\partial^2 f}{\partial z_1^2} - \frac{\lambda_1^2}{4} f = 0.\eea

The solution of the equation (\ref{xu6}) is of the form
\[
f(z_1, z_2) = c_1(z_2)e^{\frac{\lambda_1}{2}z_1}+ c_2(z_2)e^{-\frac{\lambda_1}{2}z_1}.
\]

Now putting the value of $f$ in (\ref{g1}), we get $8c_1^3(z_2)c_2(z)=\hat p_1e^{\gamma_1 z_2}$ and $8c_1(z_2)c_2^3(z)=\hat p_2e^{\gamma_2 z_2}$.

\smallskip
{\bf Sub-case 2.2.2.} Suppose that $\psi \not\equiv 0$. Clearly, from (\ref{g3}), $f$ and $\frac{\partial f}{\partial z_1}$ cannot have common zeros. 
Therefore, from (\ref{g13}), we conclude that if $z_0\in\mathbb{C}^2$ is a zero of $\frac{\partial f}{\partial z_1}$ of multiplicity $l$, then $z_0$ must be a zero of $4 \frac{\partial^2 f}{\partial z_1^2} + \lambda_1 \lambda_2 f$ of multiplicity at least $l$.
Consequently \( \psi \) is a non-zero entire function.

Thus, from (\ref{g13}), we observe that 
\[
\frac{\partial f}{\partial z_1} = 0 \Rightarrow 4 \frac{\partial^2 f}{\partial z_1^2} + \lambda_1 \lambda_2 f = 0.
\]
Therefore, $\psi$ is a non-zero entire function. Now proceeding in the same way as in the proof of Sub-case 2.1.2, we can easily deduce that $\psi$ is a non-zero polynomial in $\mathbb{C}^2$.

Now from (\ref{g13}), we obtain
\bea\label{g15}
&&\left(\psi - 4(\lambda_1 + \lambda_2)\right) f^3\\
&& = \frac{\partial f}{\partial z_1} \bigg[
4\hat a (\lambda_1 + \lambda_2) \frac{\partial f}{\partial z_1} \frac{\partial^2 f}{\partial z_1^2}
- 4\left(3f^2 + 3\hat a \left(\frac{\partial^2 f}{\partial z_1^2}\right)^2 +\hat a \frac{\partial f}{\partial z_1} \frac{\partial^3 f}{\partial z_1^3}\right)
- \lambda_1 \lambda_2\hat a \left(\frac{\partial f}{\partial z_1}\right)^2
\bigg].\nonumber\eea

Since $f$ and $\frac{\partial f}{\partial z_1}$ have no common zeros, it follows from (\ref{g15}) that
\begin{equation}\label{g16}
\psi = 4(\lambda_1 + \lambda_2),
\end{equation}
and so
\begin{equation}\label{g17}
4\hat a(\lambda_1+\lambda_2) \frac{\partial f}{\partial z_1} \frac{\partial^2 f}{\partial z_1^2}
- 4\left(3f^2 + 3\hat a \left(\frac{\partial^2 f}{\partial z_1^2}\right)^2 +\hat a \frac{\partial f}{\partial z_1} \frac{\partial^3 f}{\partial z_1^3}\right)
- \lambda_1 \lambda_2\hat a \left(\frac{\partial f}{\partial z_1}\right)^2 = 0.
\end{equation}

Now using (\ref{g14}), (\ref{g16}) and (\ref{g17}), we obtain
\[
\frac{\partial f}{\partial z_1} \left(\psi^2 \frac{\partial f}{\partial z_1} + \frac{3\hat a \lambda_1 \lambda_2}{2} f \right) = f^2 \left(1 + \frac{\hat a \lambda_1^2 \lambda_2^2}{16}\right),
\]
which implies 
\bea\label{xu7} \psi^2 \frac{\partial f}{\partial z_1} + \frac{3\hat a \lambda_1 \lambda_2}{2} f = 0.\eea

Note that $\psi$ is a non-zero polynomial and $f$ and $\frac{\partial f}{\partial z_1}$ have no common zeros. Then from (\ref{xu7}), we can easily conclude that  $\frac{\partial f}{\partial z_1}$ has no zeros, which contradicts the fact that $0$ is not Picard exceptional value of $\frac{\partial f}{\partial z_1}$.

Hence the proof is complete.

\vspace{0.1in}
{\bf Compliance of Ethical Standards:}\par

{\bf Conflict of Interest.} The authors declare that there is no conflict of interest regarding the publication of this paper.\par

{\bf Data availability statement.} Data sharing not applicable to this article as no data sets were generated or analysed during the current study.

\end{document}